\documentclass[a4paper,11pt]{amsart}
\usepackage{amsmath,inputenc,euscript,amssymb,geometry}
\usepackage{graphicx}
\usepackage{amssymb}
\usepackage{latexsym}
\usepackage{amssymb,amsbsy,amsmath,amsfonts,amssymb,amscd,color}
\usepackage{tikz}
\usepackage{mathrsfs}
\usepackage{enumitem}
\usepackage{hyperref}

\DeclareFontFamily{U}{matha}{\hyphenchar\font45}
\DeclareFontShape{U}{matha}{m}{n}{
  <-6> matha5 <6-7> matha6 <7-8> matha7
  <8-9> matha8 <9-10> matha9
  <10-12> matha10 <12-> matha12
  }{}
\DeclareSymbolFont{matha}{U}{matha}{m}{n}
\DeclareMathSymbol{\Lt}{3}{matha}{"CE}

\newtheorem{lemma}{Lemma}[section]
\newtheorem{theorem}[lemma]{Theorem}

\newtheorem{rem}[lemma]{Remark}

\numberwithin{equation}{section}

\begin{document}
\title{Remarks on hypoelliptic equations}
\author[Valeria Banica]{Valeria Banica}
\address[Valeria Banica]{Sorbonne Universit\'e, CNRS, Universit\'e de Paris, Laboratoire Jacques-Louis Lions (LJLL), F-75005 Paris, France}
\email{Valeria.Banica@sorbonne-universite.fr} 

\author[Nicolas Burq]{Nicolas Burq}
\address[Nicolas Burq]{Laboratoire de mathématiques d'Orsay, CNRS, Universit\'e Paris-Saclay, B\^at.~307, 91405 Orsay Cedex, France, and Institut Universitaire de France (IUF)}
\email{Nicolas.Burq@universite-paris-saclay.fr } 
\begin{abstract}
Many results of smooth hypoellipticity are available for scalar equations.
Much remains to be done for systems and/or  at different levels of regularity and in particular for $L^1$-hypoellipticity. In this article we provide some examples and counter-examples. 
\end{abstract}
\maketitle
\section{Introduction}
Many results are available for scalar equations that are hypoelliptic with respect to the $\mathcal C^\infty$-regularity. The general picture of $\mathcal C^\infty-$hypoellipticity for systems is unclear, and we shall give some examples and counter-examples.\smallskip

In what concerns the $L^1$-regularity, we recall that classical microlocal analysis theory (negative order pseudodifferential operators are bounded on $L^1$) yields immediately $L^1$-elliptic regularity results as
$$Au\in L^1,\quad A \mbox{ elliptic of order } k,\quad u\in \mathcal S'\quad \Longrightarrow \quad \langle D\rangle^s u\in L^1, \quad\forall s, s<k.$$
It is well known that the Hilbert transform, which is an operator of the set $\Psi^0$ of pseudodifferential operators of order zero, is not continuous on $L^1$ and hence the condition $s<k$ is sharp; see also other simple counter-examples \cite[Example 7.5]{GM}.
Restricting to Radon measures $\mathcal M$, the link between possible directions of the singular part of the measure constrained by an operator to end up at $L^1$-regularity and the non-elliptic directions of the operator were first understood by De Philippis and Rindler in \cite{DePRi} by using Murat-Tartar's elliptic wave cone. Then the following critical $L^1-$elliptic regularity result was obtained in \cite{BBWFL1}, by using the approach in \cite{DePRi} reinforced by microlocal analysis techniques and some extra geometric measure theory arguments:
\begin{equation}\label{psi0L1}
\Psi^0L^1\cap\mathcal M_0\subset L^1,
\end{equation}
where $\mathcal M_0$ denotes the finite Radon measures set, yielding for instance\footnote{We can also relax the finiteness hypothesis on $\mu$ to $\mu\in\mathcal M$ and the constraint $\Psi^0 L^1$ to $\Psi^0 (L^1_{loc}\cap \mathcal S')$ and obtain a local result:
$$\Psi^0(L^1_{loc}\cap \mathcal S')\cap\mathcal M\subset L^1_{loc}.$$
Similarly we also have 
$$\Delta\mu\in L^1_{loc}\cap\mathcal S',\quad \mu\in \mathcal S',\quad D^2\mu\in \mathcal M\quad\Longrightarrow\quad D^2\mu\in L^1_{loc}.$$}

$$\Delta\mu\in L^1,\quad \mu\in \mathcal S',\quad D^2\mu\in \mathcal M_0\quad\Longrightarrow\quad D^2\mu\in L^1.$$
In the case of hypoelliptic equations, the results in \cite{DePRi} and \cite{BBWFL1} give no satisfactory information. The problem of understanding the link between the possible directions of the singular part of the measure constrained by an operator to end up at $L^1$-regularity and the non-hypoelliptic directions of the operator was treated recently in \cite{A-sandwich} in the context of Carnot groups.  
Here we show that the microlocal analysis point of view can be extended to hypoelliptic scalar equations, to obtain higher $L^1$-regularity results. Also, we shall give some explicit examples in the vectorial case that should enlighten a conjecture on the hypoelliptic wave cone made privately to the authors by De Philippis and Rindler~\cite{DePRiprep} (see Remark \ref{conj}).  

We start with the scalar case and consider for simplicity Grushin's operator
$$G=\partial_{x}^2+x^2\partial_{y}^2,$$ 
a simple classical example of non-elliptic operator satisfying H\"ormander's hypoellipticity condition at second order, with one derivative gain that can be obtained easily: $Gu\in L^2\Rightarrow u\in H^1$ . 
We denote by $g$ the symbol of the operator $G$:
$$g(x,y,\xi,\eta)=-\xi^2-x^2\eta^2.$$
We suppose $\mu\in \mathcal M_0(\mathbb{R}^2,\mathbb R)$ solution of the following equation:
\begin{equation}\label{Gru}
G\mu=f,
\end{equation}
with $f\in L^1$.
Radon-Nikodym theorem allows for the decomposition $d\mu=\frac{d\mu}{d|\mu|}d|\mu|=hd\mathcal L^2+\frac{d\mu}{d|\mu|}d|\mu|_s$, where the non negative measure $|\mu|$ is the total variation of the measure $\mu$, the function $\frac{d\mu}{d|\mu|}\in L^1_{loc}(\mathbb{R}^2 (d|\mu|), \{\pm 1\})$ is the polar function of $\mu$, $h\in L^1(\mathbb R^2,\mathbb R)$  and the measures $|\mu|_s$ and $\mathcal L^2$ are mutually singular.
Then, see \cite{DePRi}, or simply use H\"ormander's theorem if $f\in\mathcal C^\infty$ (which gives $\mu\in C^\infty$), we have the following information on the singular part of $\mu$

$$\text{supp } |\mu|_s\overset{|\mu|_s-\text{sure}}\subseteq \{(x,y)\in\mathbb R^2, \exists (\xi,\eta)\in\mathbb S^1, g(x,y,\xi,\eta)=0\}= \{0\}\times\mathbb R,$$
yielding as an information that $\mu$ is $L^1_{loc}$ outside the line $\{(0,y),y\in\mathbb R\}$. 
First we note that by Sobolev embeddings and simple computations one can actually get $ \mu\in L^1_{loc}$. 
 Furthermore we note that even in this non-elliptic setting we can improve the conclusion to $L^1$-higher-regularity by strenghthening the microlocal analysis approach by some classical results, as follows. 

\begin{theorem}\label{th}Let  $\mu\in \mathcal {S}'(\mathbb{R}^2,\mathbb R)$ solution of \eqref{Gru} with $f\in\mathcal M_0(\mathbb R^2,\mathbb R)$. Then for all $s<1$ we have:
$$\langle D\rangle^{s}\mu\in L^1_{loc}(\mathbb R^2,\mathbb R).$$
\end{theorem}
\begin{rem}
It might be possible to get
$\mu\in W^{1,1}_{loc}$, for solutions of \eqref{Gru} with $f\in L^1$, but with a much more complicated analysis, using geometric measure theory arguments and microlocal refinements as in \cite{BBWFL1}. 
Also, we made a {\em global} assumption in~\eqref{Gru} ($f \in \mathcal{M}_0$), which at the price of little complication could be weakened ($f\in \mathcal{M}$).

\end{rem}

Now we turn our attention to non-elliptic systems, and start with a simple operator of order one, having $G$ as determinant:
$$A=\left(\begin{array}{cc}\partial_{x} & x\partial_{y}\\ -x\partial_{y} & \partial_{x}\end{array}\right).$$ 
In what concerns the $\mathcal C^\infty$-regularity, we have the following result.
\begin{theorem}\label{thvectcinfty}
i) Let  $u\in \mathcal {S}'(\mathbb{R}^2,\mathbb R^2)$ constant polarized $u=\lambda v,\lambda\in \mathbb S^1, v\in \mathcal S'(\mathbb{R}^2,\mathbb R)$. Then
$$Au\in L^2_{loc} \Longrightarrow u\in L^2_{loc},$$
and the result is sharp in the sense that we do not gain derivatives in $L^2$. Therefore $Au\in\mathcal C^\infty \Longrightarrow u\in\mathcal C^\infty,$ so $A\lambda$ is $\mathcal C^\infty-$hypoelliptic (with zero derivative gain).\medskip

\noindent
ii) There exist solutions $u\in \mathcal S'(\mathbb{R}^2,\mathbb R^2)$  to $Au=0_{\mathbb R^2}$ that are not constant polarized and do not belong to $\mathcal C^\infty$, so the operator $A$ is not $\mathcal C^\infty-$hypoelliptic.
\end{theorem}

In what concerns the $L^1$-regularity, we consider $\mu:=(\mu_1,\mu_2)\in \mathcal M_0(\mathbb{R}^2,\mathbb R^2)$ solution of
\begin{equation}\label{Gruvect}
A\left(\begin{array}{c}\mu_1\\ \mu_2\end{array}\right)=\left(\begin{array}{c}f_1\\ f_2\end{array}\right),
\end{equation}
with $f_1,f_2\in L^1$. 
From \cite{DePRi} we have the following information on the singular part of $\mu$
$$
(x,\frac{d\mu}{d|\mu|}(x))\overset{|\mu|_s-a.e}{\in}\{(x,y,\lambda_1,\lambda_2)\in\mathbb R^2\times \mathbb S^1, \exists (\xi,\eta)\in\mathbb S^1, $$
$$\left(\begin{array}{cc}\xi & x\eta\\ -x\eta & \xi\end{array}\right)\left(\begin{array}{c}\lambda_1\\ \lambda_2\end{array}\right)=0_{\mathbb R^2}\}=\{(0,y,\lambda_1,\lambda_2)\in\mathbb R^2\times \mathbb S^1\}.
$$
Therefore from these previous results we know only that the measure $\mu$ is $L^1_{loc}$ outside the line $\{(0,x_2), x_2\in\mathbb R\}$. By using the critical $L^1$-elliptic regularity results in \cite{BBWFL1}, as well as microlocal arguments involving $L^1$-regularity, in particular a F. and M. Riesz's type of theorem due to Brummelhuis \cite{Br}, we give a complete answer on the structure of singularities of the measure $\mu$. 

\begin{theorem}\label{thvectL1} i) Let  $\mu\in \mathcal M_0(\mathbb{R}^2,\mathbb R^2)$ constant polarized $\mu=(\lambda_1,\lambda_2)\nu, \lambda\in\mathbb S^1,\nu\in \mathcal M(\mathbb{R}^2,\mathbb R)$, solution of \eqref{Gruvect} with $f_1,f_2\in L^1$. Then
$$\mu\in L^1(\mathbb R^2,\mathbb R^2),$$
i.e. $A\lambda$ is $L^1-$hypoelliptic for all $\lambda\in\mathbb S^1$. 
\medskip

\noindent
ii) The result holds without the constant polarization hypothesis, also for \eqref{Gruvect} with $\mu,f_1,f_2 \in \mathcal M(\mathbb{R}^2,\mathbb R^2)$, and moreover we obtain for all $s<\frac 12$ if $\langle D\rangle^{s}\mu\in \mathcal M(\mathbb{R}^2,\mathbb R^2)$:
$$\langle D\rangle^{s}\mu\in L^1_{loc}(\mathbb R^2,\mathbb R^2).$$
\end{theorem}
\begin{rem}\label{conj}
Theorem \ref{thvectL1} ii) shows that the conjecture of G. De Philippis and F. Rindler:
\begin{equation}\label{conjDePRi}P\mu=0\Longrightarrow \frac{d\mu}{d|\mu|}(x)\overset{|\mu|_s-a.e}{\in}\{\lambda\in\mathbb S^1, P\lambda \mbox{ is not hypoelliptic}\},\end{equation}
is satisfied for $P=A$, and in this case Theorem \ref{thvectL1} i) provides an example where the right hand side set is computed.
\end{rem}

\begin{rem}
Theorem \ref{thvectL1} i) can be generalized to more general systems, that reduce to scalar equations $X_j \mu\in L^1_{loc}$ with $X_j$ involving the polarisation values $\lambda_k$, provided that the spaces spanned by the increasing brackets of the operators $X_j$'s satisfy a saturation property. Quite likely Theorem \ref{thvectL1} ii) can be also extended to more general systems. 
\end{rem}

\begin{rem}
The question of the "shape of singularities" of constant polarization measures constrained to vanish under the action of the differential operator with constant coefficients was considered in \cite{DePRi17}. It was shown that the singular part is invariant under directions orthogonal to the characteristic set of the system, if the system is of order one. This is due to the fact that the study reduces to a system of transport equations on the scalar measure $\nu$.
In \cite{BBWFL1} non-constant polarization measures subject to a system with diagonal part made by smooth vector fields were proved to have singular part and polarization propagating in a way related to the bicharacteristic flow of the system. All these results are on invariance of the set of singularities but do not give information when singularities do not exist. Theorem \ref{thvectL1} is a result in this direction, for a simple variable coefficient non-elliptic system of order one. 
\end{rem}
\begin{rem} Theorems~\ref{thvectcinfty} ii) and~\ref{thvectL1} ii) show that there exists sytems which are $L^1$ but not $C^\infty$ hypoelliptic. It is a natural and interesting question to ask whether there exist (smooth) {\em scalar} operators which have the same property.
\end{rem}
To complete the picture of various behaviors, we notice in \S 3.3 that the system 
$$\left(\begin{array}{cc}\partial_{x} & \partial_{y}\\ -x^2\partial_{y} & \partial_{x}\end{array}\right)$$ 
is both $\mathcal C^\infty$ and $L^1-$hypoelliptic.

{\bf{Aknowledgements: }} Both authors are grateful to the Institut Universitaire de France for the ideal research conditions offered by their memberships. The first author was also partially supported by the French ANR project SingFlows ANR-18-CE40-0027, while the second author was also partially supported by French ANR project ISDEEC  ANR-16-CE40-0013. The authors would like to thank the anonymous referee of \cite{BBWFL1} for pointing out the fact that the $L^1$-hypoelliptic framework is a challenging question in the field. We are also grateful to a referee of the present paper for drawing our attention to the reference \cite{NS}, giving an improvement in Theorem~\ref{th}.

\section{Proof of Theorem \ref{th}}


 Grushin operator $G$ is known to be invertible with inverse $P$ belonging to the Nagel-Stein's class $S ^{-2}_\rho$ (\S 7 d) in \cite{NS}). Applying $P$ to \eqref{Gru} we get
 \begin{equation}\label{Gruinv}
\mu=P f.
\end{equation}
For $1<p<\infty$ the class $S ^{-2}_\rho$ sends $L^p(\mathbb R^2)$ into $W^{1,p}(\mathbb R^2)$ (Theorem 1 in \cite{NS}). By duality $P$ sends $W^{-1,\overline p}(\mathbb R^2)$ into $L^{\overline{p}}(\mathbb R^2)$. Thus for $0\leq \epsilon\leq 1$ by interpolation $P$ sends $W^{-\epsilon,p_\epsilon}(\mathbb R^2)$ into $W^{1-\epsilon,\overline{p_\epsilon}}(\mathbb R^2)$ where $p_\epsilon=\frac{p}{1-\epsilon(2-p)}$ for $1<p<2$. On the other hand $L^1(\mathbb R^2) \hookrightarrow W^{-\delta,q}(\mathbb R^2)$ for all $1<q<\frac 2{2-\delta}$. It follows that one can choose $\epsilon\in[0,1]$ and $p\in [1,2]$ such that $1<p_\epsilon<\frac 2{2-\epsilon^-}$ so that we have
$$ \mathcal M_0(\mathbb R^2) \hookrightarrow W^{-0^+,1} (\mathbb R^2)\overset{Sob}{\hookrightarrow}W^{-\epsilon,p_\epsilon}(\mathbb R^2)\overset{P}{\rightarrow}W^{1-\epsilon, \overline{p_\epsilon}}(\mathbb R^2) \overset{\mathcal C^\infty_0}{\rightarrow}W^{1-\epsilon,1}(\mathbb R^2),$$

Therefore for all $s<1$ we have $\langle D\rangle^s \mu\in L^1_{loc}(\mathbb R^2)$. 

\begin{rem} The property we used  from~\cite{NS} seems to be quite specific to Grushin operator. A more general (but less precise) result can be obtained by noticing that Grushin operator $G$ satisfies H\"ormander's hypoellipticity condition at second order, and one can use the results in \cite{CCX} that yield an inversion operator $B$ in the exotic class of pseudodifferential operators $\Psi^{-1}_{\frac 12,\frac 12}$ and the reminder operator $R$ is infinitely regularizing. 
For $s<\frac 12$, since $\langle D\rangle^{s}B\langle D\rangle^\epsilon \in \Psi^{s-1+\epsilon}_{\frac 12,\frac 12}(\mathbb R^2)$ (Theorem 8.4.3 combined with the remark at the end of \S VIII 8.4 in \cite{Ho}\footnote{Note that this is not the case for Nagel-Stein's class since $\langle D\rangle ^{0^+}$ does not belong to $S^{0^+}_\rho$.}) and $\Psi^m_{\frac 12,\frac 12}$ is bounded on $L^p(\mathbb R^2)$ for $-\frac 12<m\leq 0, |\frac 1p-\frac 12|\leq |m|$ (\S VII 5.12 e), g) in \cite{Stein}), for a choice $\epsilon\in]0,1-2s]$ we have:
$$ \mathcal M_0(\mathbb R^2) \overset{\langle D\rangle^{-\epsilon} }{\rightarrow} W^{\epsilon^-,1} (\mathbb R^2)\overset{Sob}{\hookrightarrow}L^{\frac{2}{2-\epsilon^-}}(\mathbb R^2)\overset{\Psi^{s-1+\epsilon}_{\frac 12,\frac 12}}{\rightarrow}L^{\frac{2}{2-\epsilon^-}}(\mathbb R^2) \overset{\mathcal C^\infty_0}{\rightarrow}L^1(\mathbb R^2),$$
so $ \langle D\rangle^{s} B\langle D\rangle^\epsilon \langle D\rangle ^{-\epsilon} f\in L^1_{loc}(\mathbb R^2)$. However using these sharp properties of the general class $\Psi^m_{\frac 12,\frac 12}$ allows to gain only $\frac {1}{2}^-$ derivatives, while in the particular case of Grushin's inverse, Nagel-Stein's class framework ensures a better behavior. 
\end{rem}

\section{Proof of Theorems~\protect{\ref{thvectcinfty}} and \protect{\ref{thvectL1}}}\label{secthyposyst}
\subsection{Proof of Theorem \ref{thvectcinfty} i) and Theorem \ref{thvectL1} i)}\label{sectpol}
If $u=(\lambda_1,\lambda_2)v$ is a solution of 
$$\left(\begin{array}{cc}\partial_{x} & x\partial_{y}\\ -x\partial_{y} & \partial_{x}\end{array}\right)\left(\begin{array}{c}\lambda_1\\ \lambda_2\end{array}\right)v=\left(\begin{array}{c}f_1\\ f_2\end{array}\right),$$
as $(\lambda_1,\lambda_2)\in\mathbb S^1$ we get
$$\left\{\begin{array}{c}\partial_{x}v=\lambda_1 f_1+\lambda_2 f_2,\\ x\partial_{y}v=\lambda_2 f_1-\lambda_1 f_2.\end{array}\right.$$
Then we use
$$\partial_{y}=[\partial_{x},x\partial_{y}],$$
to obtain
$$\partial_{y}v=\partial_{x}(\lambda_2 f_1-\lambda_1 f_2)-x\partial_{y}(\lambda_1 f_1+\lambda_2 f_2).$$
Finally we compute
$$\Delta v=\partial_{x}(\lambda_1 f_1+\lambda_2 f_2)+\partial_{xy}(\lambda_2 f_1-\lambda_1 f_2)-x\partial_{y}^2(\lambda_1 f_1+\lambda_2 f_2).$$

At the $L^2$ level, if $f_1,f_2\in L^2_{loc}$ the inversion gives us that locally $v$ belongs to $\Psi^0 L^2\subset L^2$ (and no more information on the derivatives of $v$). Iterating the process we obtain the hypoellipticity property, as $\Delta$ commutes with $\partial_{x}$ and $\partial_{y}$. Thus we get Theorem \ref{thvectcinfty} i).

To prove Theorem \ref{thvectL1} i) we proceed similarly and since $f_1,f_2\in L^1$,  the inversion gives
$$\nu\in\Psi^0 L^1.$$
Then we use \eqref{psi0L1} to get the conclusion of Theorem \ref{thvectL1} i): $\nu\in L^1.$ 

Note that with this method we do not recover higher derivatives in $L^1$.

\subsection{Proof of Theorem \protect{\ref{thvectcinfty}} ii)}
We shall show that $A$ is not $\mathcal C^\infty$-hypoelliptic.  In this particular case, we will proceed  by an explicit computation (see below). However notice that it is actually part of a very classical line of research and our explicit calculations can be explained  by the classical caracterization of solvable operators from H\"ormander~\cite{Ho5} and the  duality  argument ``$P$ hypoelliptic iff $P^*$ solvable". We refer also to Hans Lewy's counter-example~\cite{HL} and to \cite{NiTr}. 

Consider $u_1,u_2$ defined in Fourier with respect to the second variable by
$$\hat{u_1}(x,\eta)=\chi(\eta)e^{-\frac{x^2}2\eta},\quad \hat{u_2}=-i\hat{u_1},$$
where $\chi$ has support in $\eta>0$. Then
$$\partial_{x}\hat{u_1}(x,\eta)+ix_1\eta\hat{u_2}(x,\eta)=\hat{u_1}(x,\eta)(-x\eta+i(-i)x\eta)=0,$$
$$\partial_{x}\hat{u_2}(x,\eta)-ix_1\eta\hat{u_1}(x,\eta)=\hat{u_1}(x,\eta)((-i)(-x\eta)-ix\eta)=0,$$
so $u=(u_1,u_2)$ is a solution of $Au=0_{\mathbb R^2}$. As
$$u_1(x,y)=\int e^{iy\eta} e^{-\frac{x^2}2\eta}\chi(\eta)d\eta,$$
then if $u_1\in L^2(\mathbb R^2)$ we have
$$\|u_1\|_{L^2}^2=\int\int |\hat{u_1}(x,\eta)|^2dxd\eta=\int\int e^{-y^2}|\chi(\eta)|^2dy\frac{d\eta}{\sqrt{\eta}}=C\int |\chi(\eta)|^2\frac{d\eta}{\sqrt{\eta}}.$$
Therefore by choosing $\chi$ accordingly we can obtain solutions $u_1\notin L^2$. Moreover, since
$$u_1(0,y)=\int e^{iy\eta}\chi(\eta)d\eta,$$
by taking $\chi=I_{\eta>0}$ we get Theorem \ref{thvectcinfty} ii) with a solution having as a very singular part:
$$(\Re u_1,-i\Im u_1)(0,y)=(C\delta_0(y),\frac{\tilde C}{y}).$$
We note that $(u_1e^{i\theta},-iu_1e^{i\theta})$ for $\theta\in\mathbb R$ is also a solution of $Au=0_{\mathbb R^2}$ that gives a solution for Theorem \ref{thvectcinfty} ii) with polarization on the singular part rotated by $\theta$. Therefore we get solutions for Theorem \ref{thvectcinfty} ii) with any direction of polarization on the singular part. And also in the construction $\chi=I_{\eta>0}$ can be replaced by $\chi=I_{\eta>y_0}$ for any $y_0\in\mathbb R^+$.

\subsection{Proof of Theorem \protect{\ref{thvectL1}} ii)}
We first note that the previous counterexample involves the principal value distribution of order one, which is not a measure. 

Restricting to measures we will eventually fall in the framework of F. and M. Riesz's type of theorems. The original one in the periodic setting $\mathbb T(dx)$ states that if $\mu$ is a measure and $\mbox{supp }\hat{\mu}\subseteq \mathbb N$, then $\mu$ and $dx$ have the same null sets, and in particular $\mu\Lt dx$. The equivalent on the line is that if $\mbox{supp }\hat{\mu}\subseteq (0,\infty) $ then $\mu\Lt \mathcal L^1$. In higher dimensions we have, if  $\mbox{supp }\hat{\mu}\subseteq \{\xi,\xi\cdot\omega>0\}$ for some $\omega\in\mathbb S^{d-1}$, that $\mu$ is quasi-invariant in the direction of $\omega$, thus $\mu\Lt\mathcal L^d$ (\cite{F}). 
So if the restriction of the previous type of counterexample to $x_1=0$ is a measure, then it is regular. One can have in mind also the following F. and M. Riesz's type of theorems (\cite{Sh}):
\begin{itemize}
\item if $\mu\in\mathcal M_0(\mathbb R)$ and $\hat\mu\in L^2((-\infty,0])$ then $\mu\Lt\mathcal L^1$,
\item  if $\omega\in\mathbb S^{d-1}$ and $\mu\in\mathcal M(\mathbb R^d)$ is compactly supported and $\hat\mu\in L^2(\{\xi,\xi\cdot\omega< 0\})$ then $\mu\Lt\mathcal L^d$.
\end{itemize}
In here we rely on a microlocal F. and M. Riesz's type of theorems due to Brummelhuis \cite{Br}:
\begin{theorem}  Consider $\mu\in\mathcal M(\mathbb R^d,\mathbb C)$ such that 
\begin{equation}\label{Br}
 WF_z(\mu)\cap (-WF_z(\mu))=\varnothing,\quad \forall z\in\mathbb R^d.
\end{equation}
Then $\mu\in L^1_{loc}(\mathbb R^d,\mathbb C)$.
\end{theorem}
We start first with the short proof in the case $A\mu\in\mathcal C^\infty$, and then give the more micro-local involved proof of the case $A\mu\in \mathcal M$.

\subsubsection{The case \protect{$A\mu\in\mathcal C^\infty$}}
From \eqref{Gruvect} we get the following equation on $\nu=\mu_1+i\mu_2$:
$$P\nu=F\in \mathcal C^\infty,\quad P:=\partial_x-ix\partial_y.$$
H\"ormander's theorem implies
$$WF(\nu)\subseteq Char P=\{(x,y,\xi,\eta),p(x,y,\xi,\eta)=0\}=\{(0,y,0,\eta),y,\eta\in\mathbb R\}.$$
The localization of the wave front set can be refined by removing the hypoelliptic set,\footnote{The idea is the following: as $$\|P u\|_{L^2}^2=\|\Re P u\|_{L^2}^2+\|\Im P  u\|_{L^2}^2+\langle i[\Re P,\Im P] u, u\rangle,$$
and as the principal symbol of the operator of order one $i[\Re P,\Im P]$ is $\{\Re p,\Im p\}$, one obtains a $H^\frac 12$-control of the microlocalization of $u$ on $Hyp(P)$ thus by iterations $\mathcal C^\infty$-regularity.} see \cite[Theorem 26.3.5]{Ho4}:
$$WF(\nu)\cap Hyp P=\varnothing,$$
where
$$Hyp P=\{(x,y,\xi,\eta)\in Char P,\{\Re p,\Im p\}>0\}=\{(0,y,0,\eta), y\in\mathbb R,\eta <0\}.$$
Therefore
$$WF(\nu)\subseteq\{(0,y,0,\eta), y\in\mathbb R, \eta >0\},$$
which implies
$$WF_{(x,y)}(\nu)\cap (-WF_{(x,y)}(\nu))=\varnothing,\,\,\forall (x,y)\in\mathbb R^2,$$
so applying Brummelhuis's result \eqref{Br} we obtain $\nu\in L^1_{loc}$ and so $\mu_1,\mu_2\in L^1_{loc}$.\\

\subsubsection{The case $A\mu\in \mathcal M$} In the following, localizing in space will be harmless so we can suppose $A\mu\in \mathcal M_0$. We have
$$P\nu=F\in \mathcal M_0.$$
We shall use a partition of unity $\{\chi_j\}_{0\leq j\leq 4}$ in the phase variables satisfying:
$$supp\, \chi_0\subseteq B(0,2), supp\,\chi_j\subseteq ^cB(0,1), \forall 1\leq j\leq 4,$$
and for $1\leq j\leq 4$ the cut-offs $\chi_j$ are conical, satisfying
$$ supp\, \chi_1\subseteq \{4\xi>|\eta|\},  supp\, \chi_2\subseteq \{4\xi<-|\eta|\}, $$
$$ supp\, \chi_3\subseteq \{\eta<-2|\xi|\},  supp \,\chi_4\subseteq \{\eta>2|\xi|\}.$$
We define
$$\nu_j=\chi_j((x,y),D_{(x,y)})\nu,$$
and we have
$$\nu=\sum_{j=0}^4\nu_j, \quad P\nu_j=[P,\chi_j]\nu_j+\chi_jF.$$
The first piece $\nu_0$ is microlocalized at small frequencies thus is $\mathcal C^\infty$ and in particular $L^1_{loc}$.
The following two pieces are microlocalized far from the characteristic set, where the operator $P$ is elliptic, thus we recover $W^{1^-,1}$-regularity. 
Once we shall prove that $\nu_3\in L^1_{loc}$ then we get that $\nu_4=\nu-\sum_{j=0}^3\nu_j$ is a measure, microlocalized on $\{(x,y,\xi,\eta),\eta >0\}$, to which we apply Brummelhuis's result and get $\nu_4\in L^1_{loc}$. Therefore $\nu\in L^1_{loc}$ and the conclusion of Theorem \ref{thvectL1} ii) follows for $s=0$. 
Summarising, we are left to prove that $\nu_3\in L^1_{loc}$. As 
$$P\nu_3=[P,\chi_3]\nu_3+\chi_3F,$$
by applying $\langle D_y\rangle^{-\delta}$, that belongs to $\Psi^{-\delta}$ when applied to distributions microlocalized on $|\eta|>2|\xi|$, for $0<\delta$ to be chosen later less than $\frac 12$, yields 
$$\langle D_y\rangle^{-\delta} P\nu_3=\tilde F\in\Psi^{-\delta}\mathcal M_0\subseteq W^{-\delta^+,1}\subseteq L^1$$
We note that by considering the Fourier transform in the second variable (and localizing at $|x|<1$ to simplify, as $Char P$ concerns only $x=0$)
$$supp\, \hat {\nu_3}(x,\eta),\,supp\, \hat {\tilde F}(x,\eta)\subseteq \{|x|<1,\eta<-1\},$$
and the equation on $\nu_3$ writes
$$(\partial_x+x\eta)\hat\nu_3(x,\eta)=\eta^\delta\hat {\tilde F}(x,\eta).$$
We consider now, for $p\geq 2$,
$$\nu_3^p:=\chi_p(D_y)\nu_3,$$
with $\chi_p\in\mathcal C^\infty$ satisfying  $\chi_p(\eta)=0$ for $\eta<-p-1$ and $\eta>0$, and $\chi_p(\eta)=1$ for $-p<\eta<-1$.  These distributions are compactly supported in the Fourier variables as $|\eta|<p$ implies also $|\xi|<\frac p2$. Thus $\nu_3^p\in\mathcal C^\infty$ and we have
$$\nu_3^p\rightharpoonup\nu_3.$$
We shall prove that we have $L^1$-convergence, so that in particular $\nu_3\in L^1_{loc}$. For $p>q$ we integrate the equation of $\hat{\nu_3^p}-\hat{\nu_3^q}$ 
 from $x$ to $1$ for $x>0$ (similarly we can deal with the case $x<0$ by integrating from $-1$ to $x$) to get
$$(\nu_3^p-\nu_3^q)(x,y)1_{x>0}=-\int_x^1\int e^{i(y-y')\eta+\frac {x'^2-x^2}2\eta}\eta^\delta \chi_{pq}(\eta)\tilde F(x',y')d\eta dy'dx',$$
where $\chi_{pq}=\chi_p-\chi_q$ is a localization in $\eta$ in $(-p-1,-q)$.
As $\nu_3^p,\nu_3^q$ are smooth, and localized in $|x|<1$, to get the convergence of $\nu_3^p-\nu_3^q$ in $L^1_{loc}$ it is enough to prove that the kernel
$$K_{pq}(x,y,x',y'):=1_{(x,1)}(x')\varphi(x)\int e^{i(y-y')\eta+\frac {x'^2-x^2}2\eta}\eta^\delta\chi_{pq}(\eta)d\eta,$$
where $\varphi\in\mathcal C^\infty$ is such that $\varphi(x)=1$ for $|x|<1$, $\varphi(x)=0$ for $|x|>2$,
satisfies
$$\sup_{x',y'}\|K_{pq}\|_{L^1_{x,y}}\overset{p,q\rightarrow\infty}{\longrightarrow} 0.$$
As $0<x'^2-x^2<1$ and $\eta<0$, performing $N\in\mathbb N$ integrations by parts yields
$$|K_{pq}(x,y,x',y')|\leq \frac{C}{(2|y-y'|+(x'^2-x^2))^N}\int e^{\frac {x'^2-x^2}2\eta}|\partial_\eta^N(\eta^\delta\chi_{pq}(\eta))|d\eta,$$
so for $N=0$ we get:
$$|K_{pq}(x,y,x',y')|\leq C\frac{e^{-q(x'^2-x^2)}}{(x'^2-x^2)^{1+\delta}}.$$
and for $N=1$ we get:
$$|K_{pq}(x,y,x',y')|\leq C\frac{e^{-q(x'^2-x^2)}}{(2|y-y'|+(x'^2-x^2))(x'^2-x^2)^{\delta}},$$
and for $N=2$ we get:
$$|K_{pq}(x,y,x',y')|\leq C\frac{e^{-q(x'^2-x^2)}}{(2|y-y'|+(x'^2-x^2))^2}.$$
Thus we obtain by splitting the integration in $y$ in three regions we obtain for $0<\delta<\frac 12$:
\begin{multline}
\sup_{x',y'}\|K_{pq}\|_{L^1_{x,y}}\\
\leq \sup_{0<x'< 1,y'}\int_0^{x'}\int_{\{|y-y'|<x'^2-x^2\}\cup\{x'^2-x^2<|y-y'|<(x'^2-x^2)^{\delta}\}\cup \{(x'^2-x^2)^{\delta}<|y-y'|\}}\\
\hfill |K_{pq}((x,y,x',y')|dydx
\end{multline}
$$\leq C\sup_{0<x'< 1}\int_0^{x'}\left(\frac{1}{(x'^2-x^2)^{\delta}}+\frac{\ln (x'^2-x^2)}{(x'^2-x^2)^{\delta}}+\frac{1}{(x'^2-x^2)^{\delta}}\right)e^{-q(x'^2-x^2)}dx<\frac{C}{q^{0^+}}\overset{p,q\rightarrow\infty}{\longrightarrow} 0.$$

The same proof shows that $\langle D_y\rangle^{\frac 12^-}\nu_3\in L^1_{loc}$ and so since we are in the region $|\eta|>2|\xi|$ we also have $\langle D_{x,y}\rangle^{\frac 12^-}\nu_3\in L^1_{loc}$. Thus $\langle D_{x,y}\rangle^{\frac 12^-}\nu_4=\langle D_{x,y}\rangle^{\frac 12^-}\nu-\sum_{j=0}^3\langle D_{x,y}\rangle^{\frac 12^-}\nu_j$ is a measure, microlocalized on $\{(x,y,\xi,\eta),\eta >0\}$, to which we apply again Brummelhuis's result and get $\langle D_{x,y}\rangle^{\frac 12^-}\nu_4\in L^1_{loc}$, and the conclusion of Theorem \ref{thvectL1} ii) follows for $s<\frac 12$. 

\subsection{A genuine hypoelliptic system}
For
\begin{equation}\label{hypo}\left(\begin{array}{cc}\partial_{x} & \partial_{y}\\ -x^2\partial_{y} & \partial_{x_1}\end{array}\right)\left(\begin{array}{c}u_1\\ u_2\end{array}\right)=\left(\begin{array}{c}f_1\\ f_2\end{array}\right),\end{equation}
we get
$$\left(\begin{array}{cc}G & 0\\ -x\partial_{xy} & G\end{array}\right)\left(\begin{array}{c}u_1\\ u_2\end{array}\right)=\left(\begin{array}{cc}\partial_{x} & -\partial_{y}\\ x^2\partial_{y} & \partial_{x}\end{array}\right)\left(\begin{array}{c}f_1\\ f_2\end{array}\right).$$
If $f_1,f_2\in L^2$, from the first equation we get $u_1\in L^2$ only and we cannot conclude that $u_2\in L^2$. If we suppose $f_1,f_2\in H^1$ then we get $u_1\in H^1$ and then from the second equation we get $u_2\in L^2$. Then from the first line of the initial system we get $\partial_{x_2}u_2\in L^2$ and from the second line $\partial_{x_1}u_2\in L^2$ so we also have $u_2\in H^1$. So if $f_1,f_2\in H^1$ then $u_1,u_2\in H^1$. Therefore if $f_1,f_2\in\mathcal C^\infty$ then $u_1,u_2\in\mathcal C^\infty$ (for the $L^1$ regularity we need $f_1,f_2$ to be in $W^{s,1}$ spaces).

\end{document}